\documentclass[a4paper]{amsart}  

\usepackage{amssymb} \usepackage{amscd} \usepackage{times}

\def\zbb{\mathbb{Z}}  
  
  \def\phi{\varphi}
 \def\p1{{\mathbb{P}^1_\zbb}}

\begin{document}

\title{ Lower bounds for  sup + inf and sup * inf and an Extension of Chen-Lin result in dimension 3.}
\author{Samy Skander Bahoura}

\address{Department of Mathematics, Patras University, 26500 Patras , Greece }

\email{samybahoura@yahoo.fr, bahoura@ccr.jussieu.fr}

\date{}

\maketitle

\begin{abstract}

We give two results about Harnack type inequalities. First, on compact smooth Riemannian surface without boundary, we have an estimate of the type $ \sup +\inf $. The second result concerns the solutions of prescribed scalar curvature equation on the unit ball of $ {\mathbb R}^n $ with Dirichlet condition.

Next, we give an inequality of the type $ (\sup_K u)^{2s-1} \times \inf_{\Omega} u \leq c $ for positive solutions of $ \Delta u=Vu^5 $ on $ \Omega \subset {\mathbb R}^3 $, where $ K $ is a compact set of $ \Omega $ and $ V $ is $ s-$ h\"olderian,$ s\in ]-1/2,1] $. For the case $ s=1/2 $, we prove that if $ \min_{\Omega} u>m>0 $ and the h\"olderian constant $ A $ of $ V $ is small enough ( in certain meaning), we have the uniform boundedness of the supremum of the solutions of the previous equation on any compact set of $ \Omega $.

\end{abstract}

\bigskip

\bigskip

\begin{center}  1. INTRODUCTION AND RESULTS.
\end{center}

\bigskip

We denote $ \Delta =-\nabla^j(\nabla_j) $, the geometric Laplacian.

\medskip

On compact smooth Riemann surface without boundary $ (M,g) $ we consider the following equation :

$$ \Delta u+k=V e^u, \qquad (E_1) $$

with, $ k\in {\mathbb R}_{*,+} $ and $ 0 \leq V \leq  b $ ( $ V \not \equiv 0 $ ).

\medskip

We suppose $ V $ smooth. The previous equation is of type prescribed scalar curvature. We search to know if it's possible to have a priori estimate of the type $ \sup +\inf $.

\bigskip

Note that in dimension 2, on $ {\mathbb R}^2 $, we have different results about $ \sup + \inf $ inequalities for the following equation:

$$ \Delta u = V e^u, \qquad (E_2) $$

see [B-L-S], [B-M], [C-L 2], [L 2] and [S].

\bigskip

In [S], Shafrir proved an inequality of the type $ \sup u + C \inf u < C' $ with minimal conditions on the prescribed scalar curvature. In [B-L-S], Brezis-Li-Shafrir have proved a $ \sup u + \inf u $ inequality with lipschitzian assumption on prescribed curvature. Finaly, [C-L 2] have proved the same result with h\"olderian assumption on $ V $ in the equation $ (E_2) $.

\bigskip

Here, we are interested by the minoration of this sum. We can suppose that $ Volume(M)=1 $. We obtain,

\bigskip

{\it Theorem 1. For all  $ k,b>0 $, there exists a constant  $ c=c(k,b,M,g) $ such that, for all solution of $ (E_1) $:

 $$ \dfrac{k-4\pi}{4\pi} \sup_M u + \inf_M u \geq c. $$}

We can remark that for $ k=8\pi $, we have the same result than in [B 1]. Here there is no restriction on $ k $.

\bigskip

Now we work in dimension $ n\geq 3 $, we set $ B=B_1(0) $ the unit ball of $ {\mathbb R}^n $. We try to study some properties of the solutions of the following equation:

$$ \Delta u =Vu^{N-1-\epsilon}, u>0 \,\, {\rm in } \,\, B,\,\, u=0 \,\, {\rm on} \,\, \partial B \qquad (E_3)  $$

with $ 0 \leq V (x) \leq b < +\infty $, $ 0 \leq \epsilon < 2/(n-2) $ and $ N=\dfrac{2n}{n-2} $ the critical Sobolev exponant.

\bigskip

Equation $ (E_3) $ is the prescribed scalar curvature equation, it was studied a lot. We know, after using Pohozaev identity that, there is no solution for this equation if we assume $ \epsilon =0 $ and $ V \equiv 1 $, see [P].

\bigskip

{\it Theorem 2. For all compact $ K $ of $ B $, there exists one positive constant $ c=c(n,b,K) $ such that for all solution of $ (E_3) $ :

$$ (\sup_B u )^7 \times \inf_K u \geq c. $$}

Recall that estimates like in the last theorem exist, see for example [B 1] et [B 2].

\bigskip

Now we work on $ \Omega \subset {\mathbb R}^3 $ and we consider the following equation:

$$ \Delta u=V u^5, \,\,\, u > 0, \qquad (E_4)$$

with,

$$ 0 < a \leq V(x)\leq b \,\,\, {\rm and} \,\,\, |V(x)-V(y)|\leq A|x-y|^s, \,\,\, s\in [\dfrac{1}{2}, 1], x,y \in \Omega. \qquad (C) $$

Without loss of genarality, we suppose $ \Omega = B $ the unit ball of $ {\mathbb R}^3 $.

\bigskip

The equation $ (E_4) $ is the scalar curvature equation in three dimensions. It was studied a lot, see for example [B 3], [C-L 1], [L 1]. In [C-L 1], Chen and Lin have proved that if $ s > \dfrac{1}{2} $, then each sequence $ (u_k)_k $ which are solutions of $ (E_4) $ ( with fixed $ V $ ) are in $ L_{loc}^{\infty} $ if we suppose $ \min_B u_k >m >0 $. When $ s=1 $ they prove that the $ \sup \times \inf $ inequality holds. To prove those results, they use the moving-plane method.

\bigskip

In [L 1], Li proved (in particular) that the product $ \sup \times \inf $ is bounded if we replace $ \Omega $ by the three sphere $ {\mathbb S}_3 $. He used the notion of isolated and isolated simple blow-up points.

\bigskip

We can see in [B 3] another proof of the boundedness of $ \sup^{1/3} \times \inf $, also with the moving-plane method.

\bigskip

Note that, if we suppose $ \Omega $ a Riemannian manifold of dimension 3 (not necessarily compact), Li and Zhang (see [L-Z]) have proved that the $ \sup \times \inf $ holds when the prescribed scalar curvature is a constant.

\bigskip

Note that, in our work, we have no assumption on energy. There are many results, if we suppose the energy bounded.

\bigskip

Here, we use the moving-plane method to have $ \sup \times \inf $ inequalities. This method was developed by Gidas-Ni-Nirenberg, used by Chen-Lin and Li-Zhang, see  [G-N-N], [C-L 1] and [L-Z]. In our work we follow and use the technique of Li and Zhang, see [L-Z].

\bigskip

{\it Theorem 3. If $ s\in ]\dfrac{1}{2},1] $, \, then, \, for all positive numbers $ a,b, A $ and all compact $ K $ of $ B $, there exists a positive constant $ c=c(a, b, A, s, K) $ such that:

$$ (\sup_K u )^{2s-1} \times \inf_B u \leq c, $$

where $ u $ is solution of $ (E_4) $ with $ V $ satisfying $ (C) $.

For $ s=\dfrac{1}{2} $ and $ a,b,m >0 $, there exists $ \delta =\delta (a,b,m)>0 $ such that for $ u $ solution of $ (E_4) $ with $ A \in ]0,\delta] $ for $ V $ in $ (C) $ and $ u \geq m $, we have:

$$ \sup_K u \leq c = c(a, b, m, K ),$$

where $ K \subset \subset B_1 $.}

\bigskip

Note that in [B 3], for the dimension 4, we have a result like in the second part of the  theorem 3.

\bigskip

About usual Harnack inequalities, we can find in [G-T] lots of those estimates. For harmonic functions ($ \Delta u=-\sum_{i=1}^n \partial_{ii} u= 0 $ on open set of $ {\mathbb R}^n $), we have an estimate of the type:

 $$ \dfrac{\sup_{B_R} u}{\inf_{B_R} u} \leq 3^n, $$

on small ball $ B_R $ of radius $ R $ (see chapter 1 in [G-T]).

\bigskip

We have other results if we consider a general elliptic operator ($ L= \partial_i(a^{ij}\partial_j)+\sum_{j=1}^n b^j\partial_j +c $ on open set of $ {\mathbb R}^n $), we obtain:

$$ \dfrac{\sup_{B_R} u}{\inf_{B_R} u} \leq C[n,R,c, (b^j)_j,(a^{ij})_{i,j}] $$

for a non negative function $ u $ such that $ Lu=0 $ ( $ B_R $ is a ball of radius $ R $ . See for example theorem 8.20 in [G-T].

\bigskip

For subharmonic and superharmonic functions there is another type of Harnack inequalities linking their norm $ L^p $ to their infimum or supremum. (See chapter 8 in [G-T]).

\bigskip

Here we follow the same idea and we try to compare the $ \sup $ and the $ \inf $ in a certain meaning.

\bigskip

\bigskip

\begin{center}  {2. PROOFS OF THE THEOREMS.}
\end{center}

\bigskip

\bigskip

\underbar {\bf Proof of Theorem 1:}

\bigskip

Consider the equation :

$$ \Delta u_i+k=V_ie^{u_i}, $$

\underbar {\bf Case 1: $ \sup_M u_i \leq c < +\infty $.}

\bigskip

We set $ x_i $ the point where $ u_i $ is maximum, $ u_i(x_i)=\sup_M u_i $, then:

$$ 0 \leq \Delta u_i(x_i)=V_i(x_i)e^{u_i(x_i)}-k\leq be^{u_i(x_i)}-k, $$

thus,

$$  \log \left (\dfrac{k}{b} \right ) \leq u_i(x_i) \leq c'. $$

We denote $ G $ the Green function of laplacian,

$$ \Delta_{y, distribution } G(x,.)=1-\delta_x \,\, {\rm and} \,\, G(x,y) \geq 0, \int_M G(x,y)\equiv C. $$

we can write,

$$ \log \left (\dfrac{k}{b} \right ) \leq u_i(x_i)=\int_M u_idV_g-\int_M G(x_i,y)[V_i(y)e^{u_i(y)}-k] dV_g \leq \int_M u_i+C(be^c-k). $$

We deduce:

$$ -\infty < c_2 \leq \int_M u_i \leq c_1 < +\infty, \,\, \forall \,\, i. $$

Now, we write,

$$ \min_M u_i=u_i(y_i)=\int_M u_i+\int_M G(y_i,y)[V_i(y)e^{u_i(y)}-k]dV_g \geq c_2-kC >-\infty .$$

Thus,

$$ ||u_i||_{L^{\infty}} \leq c' <+\infty, \,\, \forall \,\, i. $$

\underbar {\bf Case 2: $ \sup_M u_i \to +\infty .$}

\bigskip

According to T. Aubin (see [A]), we have,

$$ G(x,y)=-\dfrac{1}{2 \pi} \log d(x,y)+ g(x,y), $$

where , $ g $ is a regular part of $ G $, it is a continuous function on $ M\times M $.

\bigskip

Let us note $ x_i $ the point where $ u_i $ is maximum, $ u_i(x_i)=\max_M u_i $. We can suppose that $ x_i \to x_0 $ and in the conformal isothermal coordinates around $ x_0 $ we set $ v_i(x)=u_i(x_i+xe^{-u_i ( x_i ) /2})-u_i(x_i) $, then,

$$ \Delta v_i+ h_i=\tilde V_i e^{v_i},\,\, h_i \to 0 $$

$ v_i(0)=0, v_i(x) \leq 0, 0 \leq \tilde V_i(x)\leq b $. We can use theorem 3 of [B-M] and we deduce after passing to the subsequence that:

$$ v_i(x) \geq C >-\infty, \,\,  {\rm for } \,\, |x| \leq r. $$

Thus,

$$ u_i(y) \geq u_i(x_i)+ C, \,\, {\rm if } \,\, d(y, x_i)\leq r e^{-u_i(x_i)/2}, $$

Now, we work on $ M-B(x_i, re^{-u_i(x_i)/2}) $,

$$ G(x_i,y)\leq \dfrac{1}{4 \pi} u_i(x_i)+ C_1,\,\, {\rm on} \,\, \partial B(x_i, re^{-u_i(x_i)/2})  $$

$$ \Delta [u_i-kG(x_i,.)] \geq 0, \,\, {\rm on} \,\, M-B(x_i, re^{-u_i(x_i)/2}) $$

$$ u_i(y)-kG(x_i,y)- \dfrac{ 4\pi-k}{4\pi} u_i(x_i) + kC_1 - C \geq 0, \,\, {\rm on} \,\, \partial B(x_i, re^{-u_i(x_i)/2}). $$

By maximum principle, we obtain:

$$ u_i \geq kG(x_i,.) + \dfrac{ 4\pi-k}{4\pi} u_i(x_i)-kC_1 + C, \,\, {\rm on } \,\, M-B(x_i, re^{-u_i(x_i)/2}), $$

We use the fact, $ \int_M G(x_i,y)\equiv constant $, and by integration of the last inequality we have,

$$ \inf_M u_i+ \dfrac{ k-4\pi}{4\pi} \sup_M u_i \geq c > -\infty , $$

{\bf Example with $ V_i \to 0 $ :} we can take, $ u_i \equiv \log k+\log i $ and $ V_i \equiv 1/i $.

\bigskip

{\bf Remark:} If we suppose $ V_i \geq a> 0 $ uniformly, then, when $ k < 4 \pi $ we can not have $ \sup_M u_i \to +\infty $. To see this, it is sufficient to integrate the equation.

\bigskip

\underbar {\bf Proof of Theorem 2:}

\bigskip

We are going to prove that each sequence has a subsequence who has the searched inequality.

\bigskip

Next, we use the fact that, if we have possibility to extract a subsequence we do it and we denote $ (u_i)_i $ the subsequence.

\bigskip

We have,

$$ \Delta u_i=V_iu_i^{N-1-\epsilon_i}, \,\, u_i>0 \,\, {\rm on } \,\, B, \qquad (\tilde E) $$

with $ 0 \leq V_i(x) \leq b $ ( $ V_i \not \equiv 0 $).

\bigskip

Let us note $ G $ the Green function of the laplacian on unit ball with Dirichlet condition. G is of the form:

$$ G(x,y)=\dfrac{1}{n(n-2) \omega_n |x-y|^{n-2}}-\dfrac{1}{n(n-2)\omega_n(|x|^2|y|^2+1-2x.y)^{(n-2)/2}}. $$

Denote $ x_i $ the point where $ u_i $ is maximum. We write:

$$ u_i(x_i)=\int_B G(x_i,y)V_i(x)[u_i(y)]^{N-1-\epsilon_i} dy \leq b [u_i(x_i)]^{N-1-\epsilon_i}\int_{B} G(x_i,y) dy. $$

Consider the function  $ h(x)=|x|^2-1 $, we have:

$$ \int_B G(x_i,y)dy =\dfrac{1-|x_i|^2}{2n} \leq  d(x_i,\partial B )/n. $$

We deduce:

$$ 0 < \dfrac{n}{b} \leq [u_i(x_i)]^{4/(n-2)-\epsilon_i} d(x_i, \partial B). $$

\underbar { {\bf Case 1:}} $ \max_B u_i \leq c $

\bigskip

Then, $ d(x_i, \partial B) \geq c'>0 $. By elliptic estimates, $ u_i \to u $, with $ u>0 $. Then, $ \inf_K u_i \geq \tilde c>0 $ with $ K \subset \subset B $.

\bigskip

To see this, we can write $ (\tilde E) $ as:

$$ \Delta u_i=f_i $$

with, $ f_i $ uniformly in $ L^p $ for $ p >n $. We can use the elliptic estimates to have $ u_i $ uniformly in $ W^{2,p}(B) $ and by the Sobolev embedding, we have $ u_i $ uniformly in $ C^{1,\theta}(\bar B) $, for some $ \theta \in ]0,1[ $.

\bigskip

Now, we can see that:

$$ \int_B \nabla u_i.\nabla \phi = \left ( \int_B V_i u_i^{N-1-\epsilon_i}\phi \right ) \geq 0 ,\,\, {\rm for\,\, all} \,\, \phi \in C_0^{\infty}(B), \,\, \phi \geq 0 \,\, {\rm (distribution)}. $$

We can passe to the limit $ u_i \to u \geq 0 $ (subsequence) and $ u\in C^1(\bar B) $. Then, we have:

$$ \int_B \nabla u.\nabla \phi \geq 0, \,\, \phi \in C_0^{\infty}(B), \,\, \phi \geq 0 \,\, {\rm (distribution)}. $$

We can use the strong maximum principle for weak solutions, see for example, Gilbarg-Trudinger, theorem 8.19 (applied to $ -u \leq 0 $):

\bigskip

If, there is a point $t $ in $ B $ such that, $ u (t)=0 $ then, $ u\equiv 0 $. But we can see that $ u_i(x_i) \geq \tilde c' >0 $ with $ \tilde c' $ do not depends on $ i $ and $ x_i \not \to \partial B $ (subsequence).

\bigskip

Finaly, $ u >0 $ on $ B $.

\bigskip

{\bf Remark 1}. Why do we do this ? in fact, we have neither $ u_i \in C^2(\bar B) $ nor $ u_i \to u $ in $ C^2 $ norm because we don't have more regularity on $ V_i $ and finally we don't have $ \Delta u \geq 0 $ in the strong sense. We have weakly $ \Delta u \geq 0 $ with a good regularity on $ u $. Here, it is sufficient to have:  $ C^1 $ regularity on $ u $ and an uniform boundedness for $ u_i $ in $ C^{1,\theta} $ ($ 0 < \theta <1 $), to obtain a good convergence for  $ u_i $. After we can use a strong maximum principle for weak solutions.

\bigskip

{\bf Remark 2}. If we take a sequence of functions $ V_i $ which converge uniformly to  $ 0 $ ( for example), the previous case 1 is not possible.

\bigskip

\underbar { {\bf Case 2:}} $ \max_B u_i \to + \infty $

\bigskip

\underbar { \bf { I)  $ x_i \to x_0 \in \partial B : $}}

\bigskip

To simplify our computations, we assume $ n/b > 1/2 $. Then, $ B(x_i,r_i) \in B $, with $ r_i=\dfrac{1}{2[u_i(x_i)]^{4/(n-2)-\epsilon_i}} $. We consider the following functions :

$$ v_i(x)=\dfrac{u_i[x_i+x/[u_i(x_i)]^{2/(n-2)-\epsilon_i/2}]}{u_i(x_i)}, $$

Those functions $ v_i $, exist on  $ \Omega_i=B(0,5t_i) $, $ t_i=1/10[u_i(x_i)]^{2/(n-2)-\epsilon_i/2} $. We have :

$$ \Delta v_i=\tilde V_i v_i^{N-1-\epsilon_i}, \,\, 0 < v_i(x) \leq v_i(0)=1, \,\, 0 \leq \tilde V_i(x) \leq b. $$

with, $ \tilde V_i(x)=V_i[x_i+x/[u_i(x_i)]^{2/(n-2)-\epsilon_i/2}] $.

\bigskip

We use Harnack inequality for $ v_i $ ( see Theoreme 8.20 of [G-T]), we obtain:

$$ \max_{B(0,t_i)} v_i \leq C \inf_{B(0,t_i)} v_i. $$

where $ C=[C_0(n)]^{1+b} $ ( see [G-T] and $ t_i \leq 1 $).

\bigskip

In $ 0 $, we obtain: $ u_i(x)\geq C(n,b) u_i(x_i) $ for $ |x|\leq s_i = 1/10[u_i(x_i)]^{4/(n-2)-\epsilon_i} $. Let us note that, $ C(n,b)=C=[C_0(n)]^{1+b} $.

\bigskip

If we consider $ B-B(x_i,s_i) $, then,

$$ G(x_i,y) \leq c(n) [u_i(x_i)]^{4-(n-2)\epsilon_i} , \,\, {\rm for } \,\, d(x_i,y)=s_i, $$

$$ \Delta G(x_i,.)=0, \,\, G(x_i,.)_{| \partial B} =0, $$

with, $ c(n)=\dfrac{10^{n-2}}{n(n-2)\omega_n} $.

\bigskip

Thus,

$$ u_i(y)-\dfrac{C(n,b) G(x_i,y)}{c(n) [u_i(x_i)]^{3-(n-2)\epsilon_i}} \geq 0, \,\, {\rm for} \,\, d(y,x_i)=s_i, \,\, {\rm or, \,\, on } \,\, \partial B. $$

$$ \Delta \left [ u_i-\dfrac{C(n,b) G(x_i,.)}{c(n) [u_i(x_i)]^{3-(n-2)\epsilon_i}}\right ] \geq 0. $$

By maximum principle, we have:

$$ u_i(y)-\dfrac{C(n,b) G(x_i,y)}{c(n) [u_i(x_i)]^{3-(n-2)\epsilon_i}} \geq 0, \,\, {\rm on} \,\, B-B(x_i,s_i). $$

In other terms,

$$ u_i(y) \geq \dfrac{C(n,b)}{c(n)} G(x_i,y)[u_i(x_i)]^{-3+(n-2)\epsilon_i}, \,\, {\rm on }\,\, B-B(x_i,s_i).$$

Now, we know that,

$$ G(x_i,y) \geq c'(n) (1-|y|)^{n-2}\times (1-|x_i|)^{n-2}. $$

where $ c'(n)=\dfrac{1}{2(n-2)2^{2(n-2)} \omega_n} $.

\bigskip

We denote, $ c'(n,b)=\dfrac{n}{b} $. Using the fact, $ 1-|x_i|=d(x_i,\partial B) \geq c'(n,b)[u_i(x_i)]^{-4/(n-2)+\epsilon_i} $, we obtain,

$$ u_i(y) \geq \dfrac{C(n,b)c'(n)c'(n,b)}{c(n)} (1-|y|)^{n-2} [u_i(x_i)]^{-7+2(n-2)\epsilon_i}, \,\, {\rm on} \,\, B-B(x_i,s_i). $$

On $ B(0,k) $ with $ k <1 $, by maximum principle we have: $ \inf_{B(0,k)} u_i=\inf_{\partial B(0,k)} u_i $.

\bigskip

Then,

$$  u_i(y) \geq C(n,b) (1-k)^{n-2} [u_i(x_i)]^{-7+2(n-2)\epsilon_i}, \,\, {\rm on } \,\, B(0,k)-B(x_i,s_i), $$

but, $ x_i \to x_0 \in \partial B $, and for $ i $ large we can conclude that $ B(x_i,s_i) \cap B(0,k)=\emptyset  $  and thus,

$$ \inf_{B(0,k)} u_i \times [u_i(x_i)]^{7} \geq C(n,b,k). $$

We can remark that:

$$ C(n,b,k)=\dfrac{C(n,b)c'(n)c'(n,b)}{c(n)} (1-k)^{n-2}. $$

with, $ C(n,b)=C_0(n)^{1+b} $, $ c(n)=\dfrac{10^{n-2}}{n(n-2)\omega_n} $, $ c'(n)=\dfrac{1}{2(n-2)2^{2(n-2)} \omega_n} $ and $ c'(n,b)=\dfrac{n}{b} $.

\bigskip

Then,

$$ C(n,b,k)=\dfrac{[C_0(n)]^{1+b}2n^2(n-2)2^{2(n-2)} \omega_n}{bn(n-2)\omega_n} (1-k)^{n-2}=\dfrac{[C_0(n)]^{1+b}2n2^{2(n-2)}}{b}(1-k)^{n-2}. $$

\bigskip

\underbar {\bf { II) $ x_i \to x_0 \in B $ : }}

\bigskip

Our computations are the same as in the previous case I), there are some modifications.

\bigskip

We take, $ t_i=1 $ and $ s_i= [u_i(x_i)]^{2/(n-2)-\epsilon_i/2} $. We have:

$$ G(x_i,y) \leq C(n) [u_i(x_i)]^{2-(n-2)\epsilon_i/2}. $$

After,

$$ u_i(y) \geq C'(n,b) G(x_i,y) [u_i(x_i)]^{-1+(n-2)\epsilon_i/2}, $$

we use the fact, $ x_i \to x_0 \in B $, $ G(x_i,y) \geq C''(n,b, x_0) (1-k)^{n-2}, $

\bigskip

then,

$$ \inf_{B(0,k)} u_i \times [u_i(x_i)]^{1-(n-2)\epsilon_i/2} \geq c(n,b,k,x_0)>0. $$

\bigskip

\underbar {\bf Proof of the Theorem 3}

\bigskip

\underbar {\bf Step 1: blow-up technique}

\bigskip

We are going to prove the following assertion:

$$ \exists \,\,\, c,R >0 \,\,\, {\rm such \,\, that}, \,\,\, R \left [\sup_{B(0,R)} u \right ]^{2s-1} \times \inf_B u \leq c \,\,\, {\rm if} \,\,\,  \dfrac{1}{2} < s \leq 1, $$

and,

$$ \exists \,\,\, c,R >0 \,\,\, {\rm such \,\, that}, \,\,\, R \sup_{B(0,R)} u \leq c \,\,\, {\rm if } \,\,\, s=\dfrac{1}{2} .$$

We argue by contradiction ( and after passing to a subseqence) and we suppose that for $ R_k \to 0 $ we have:

$$ R_k \left [\sup_{B(0,R_k)} u_k \right ]^{2s-1} \times \inf_B u_k \to +\infty, \,\,\, {\rm for} \,\,\, s\in ]1/2,1]. $$

$$ R_k \sup_{B(0,R_k)} u_k \to + \infty, \,\,\, {\rm for } \,\,\, s=1/2. $$

Let $ x_k $ be the point such that $ u_k(x_k)= \sup_{B(0,R_k)} u_k $ and consider the following function:

$$ s_k(x)= \sqrt { (R_k-|x-x_k|)} u_k(x). $$

Let $ a_k $ be the point such that: $ s_k(a_k)= \sup_{B(x_k,R_k)} s_k $. We set $ M_k=u_k(a_k) $ and $ l_k=R_k-|a_k-x_k| $. We have:

$$ M_k^{-1} u_k(x) \leq \sqrt 2, \,\,\, {\rm for} \,\,\, |x-a_k|\leq \dfrac{l_k}{2} M_k^2. $$

We have:

$$ \dfrac{l_k}{2}M_k \to +\infty, \,\, v_k(y)=M_k^{-1} u_k(a_k+M_k^{-2}y) \,\, {\rm for}\,\, |y| \leq \dfrac{l_k}{2}M_k^2, $$

$$ \Delta v_k=V_k{v_k}^5, \,\, v_k(0)=1, \,\, 0 < v_k \leq {\sqrt 2}. $$

We know, after passing to a subsequence, that:

$$ v_k \to U, \,\, {\rm with } \,\,  \Delta U=V(0)U^5, \,\,\, U>0, \,\,\, {\rm on } \,\, {\mathbb R}^3. $$

It easy to see that we can suppose $ V(0)=1 $. The result of Caffarelli-Gidas-Spruck (see [C-G-S] ) assures that $ U $ has an explicit form and is radially symmetric about some point.

\bigskip

\underbar { \bf Step 2: The moving plane method}

\bigskip

Now, we use the Kelvin transform and we set for $ \lambda >0 $ :

$$ v_k^{\lambda}(y)=\dfrac{\lambda}{|y|}v_k(y^{\lambda}) \,\, {\rm with } \,\, y^{\lambda}=\dfrac{\lambda^2y}{|y|^2}. $$

We denote $ \Sigma_{\lambda} $ by:

$$ \Sigma_{\lambda}=B \left (0,R_kM_k^{2s} \right)-{\bar  B(0,\lambda)}. $$

We have the following boundary condition:

$$ \lim_{k\to + \infty} \min_{|y|=R_kM_k^{2s}}(v_k(y)|y|) \to +\infty. $$

We have:

$$ \Delta v_k^{\lambda}=V_k^{\lambda}(v_k^{\lambda})^5. $$

We set:

$$ w_{\lambda}=v_k-v_k^{\lambda}. $$

Then,

$$ \Delta w_{\lambda}+\dfrac{n+2}{n-2}\xi^{4}V_kw_{\lambda}=E_{\lambda}, $$

with $ E_{\lambda}=(V_k-V_k^{\lambda})(v_k^{\lambda})^5 . $

\bigskip

Clearly, we have the following lemma.

\bigskip

\underbar {\bf Lemma 1:}

\bigskip

We have:

$$ |E_{\lambda}|\leq A_k \times C(\lambda_1)M_k^{-2s}\lambda^5 |y|^{s-5} \leq C(\lambda_1) \lambda_1^s \times A_k M_k^{-2s} \lambda^{5-s}|y|^{s-5}. $$

Let

$$ h_{\lambda}=-C(s,\lambda_1) A_k M_k^{-2s}\left [ 1-\left (\dfrac{\lambda}{|y|}\right)^{4-s} \right ]. $$

\underbar {\bf Lemma 2:}

$$ \exists \,\, \lambda^k_0>0 \,\,\, {\rm such\,\, that} \,\,\,  w_{\lambda}+h_{\lambda} > 0 \,\, {\rm in } \,\, \Sigma_{\lambda} \,\, \forall \,\, 0< \lambda \leq \lambda^k_0. $$

The proof of the lemma 2 is like the proof of the step 1 of the lemma 2 in [L-Z], we omit it here.

\bigskip

We set:

$$ \lambda^k= \sup \{ \lambda \leq \lambda_1, \,\, {\rm such \,\, that} \,\, w_{\mu}+h_{\mu} >0 \,\, {\rm in} \,\, \Sigma_{\mu} \,\, {\rm for \,\, all } \,\, 0 < \mu \leq \lambda  \}. $$

We have:

\bigskip

$  {\rm If } \,\, s\in ]1/2,1] \,\, {\rm then } \,\, |h_{\lambda^k}|R_k M_k^{2s} \leq C(s,\lambda_1) \sup_k A_k, $ and thus,

$$ w_{\lambda^k}+h_{\lambda^k} >0 \,\,\, \forall  \,\,\, |y|=R_kM_k^{2s} . $$

$  {\rm If } \,\, s=\dfrac{1}{2}, \min_M u_k \geq m>0 $ and $  A_k \to 0 $, we obtain,

$$ \min_{ |y|= (2\lambda_1 M_k)/m}[v_k(y)|y|] \geq 2\lambda_1 >0, $$

thus, for $ |y|=\dfrac{2\lambda_1 M_k}{m} $ and $ k $ large we have:

$$ w_{\lambda^k}+h_{\lambda^k} \geq \dfrac{[-\lambda^k v_k(y^{\lambda})+ 2\lambda_1-C(\lambda_1,s) A_k]}{(2\lambda_1 M_k)/m} \geq \dfrac{m}{2\lambda_1 M_k} [-(1+\epsilon )\lambda_1-\epsilon \lambda_1 + 2 \lambda_1] >0 , $$

where $ \epsilon >0 $ is very small and $ v_k(y^{\lambda}) \to U(0)= 1 $.

\bigskip

For the case $ s=\dfrac{1}{2} $, we work in $ \Sigma_{\lambda}= B \left (0,\dfrac{2\lambda_1 M_k}{m} \right )-\bar B(0,\lambda) $. It is easy to see that, $ \dfrac{2\lambda_1 M_k}{m} << R_k M_k^2 $. We define $ \lambda^k $ as in the case $ 1/2 < s\leq 1. $

\bigskip

If we use the Hopf maximum principle, we prove that $ \lambda^k=\lambda_1 $ like in [L-Z]. We have the same contradiction as in [L-Z].

\bigskip

\bigskip

\begin{center}

{\bf  ACKNOWLEDGEMENT. }

\end{center}

\smallskip

This work was done when the author was in Greece at Patras. The author is grateful to Professor Athanase Cotsiolis, the Department of Mathematics of Patras University  and the IKY Foundation for hospitality and the excellent conditions of work.

\bigskip

\bigskip

\bigskip

\underbar {\bf References:}

\bigskip

[A] T.Aubin, Some nonlinear problems in Riemannian Geometry. Springer-Verlag, 1998.

\smallskip

[B 1] S.S. Bahoura. In\'egalit\'es de Harnack pour les solutions d'equation du type courbure sclaire prescrite. C.R. Math. Acad. Sci. Paris. 341 (2005), no.1, 25-28.

\smallskip

[B 2] S.S. Bahoura. Estimations du type $ \sup u\times \inf u $ sur une vari\'et\'e compacte. Bull. Sci. Math. 1-13 (2006).

\smallskip

[B 3] S.S Bahoura. Majorations du type $ \sup u \times \inf u \leq c $ pour l'\'equation de la courbure scalaire sur un ouvert de $ {\mathbb R}^n, n\geq 3 $. J. Math. Pures. Appl.(9) 83 2004 no, 9, 1109-1150.

\smallskip

[B-L-S], H. Brezis, YY. Li and I. Shafrir. A sup+inf inequality for some
nonlinear elliptic equations involving exponential
nonlinearities. J.Funct.Anal.115 (1993) 344-358.

\smallskip

[B M] H. Brezis, F. Merle, Uniform estimates and Blow-up Behavior for Solutions of $ -\Delta u=V(x) e^u $ in two dimension. Commun. in Partial Differential Equations, 16 (8 and 9), 1223-1253(1991).

\smallskip

[C-G-S] L. Caffarelli, B. Gidas, J. Spruck. Asymptotic symmetry and local
behavior of semilinear elliptic equations with critical Sobolev
growth. Comm. Pure Appl. Math. 37 (1984) 369-402.

\smallskip

[C-L 1] C-C.Chen, C-S. Lin. Estimates of the conformal scalar curvature
equation via the method of moving planes. Comm. Pure
Appl. Math. L(1997) 0971-1017.

\smallskip

[C-L 2] C-C.Chen, C-S. Lin. A sharp sup+inf inequality for a nonlinear elliptic equation in ${\mathbb R}^2$.
Commun. Anal. Geom. 6, No.1, 1-19 (1998).

\smallskip

[G-T] D. Gilbarg, N.S. Trudinger. Elliptic Partial Differential Equations of Second order, Berlin Springer-Verlag, Second edition, Grundlehern Math. Wiss.,224, 1983.

\smallskip

[L 1] YY. Li. Prescribing scalar curvature on $ {\mathbb S}_n $ and related
Problems. C.R. Acad. Sci. Paris 317 (1993) 159-164. Part
I: J. Differ. Equations 120 (1995) 319-410. Part II: Existence and
compactness. Comm. Pure Appl.Math.49 (1996) 541-597.

\smallskip

[L 2] YY.Li, Harnack type Inequality, the Methode of Moving Planes. Commun. Math. Phys. 200 421-444.(1999).

\smallskip

[L-Z] YY. Li, L. Zhang. A Harnack type inequality for the Yamabe equation in low dimensions.  Calc. Var. Partial Differential Equations  20  (2004),  no. 2, 133--151

\smallskip

[P] S. Pohozaev, Eigenfunctions of the equation $\Delta
 u+\lambda f(u)=0 $. Soviet. Math. Dokl., vol. 6 (1965), 1408-1411.

\smallskip

[S] I. Shafrir. A sup+inf inequality for the equation $ -\Delta u=Ve^u $. C. R. Acad.Sci. Paris S\'er. I Math. 315 (1992), no. 2, 159-164.

\end{document}